\documentclass[a4paper,12pt]{article}
\usepackage{amsmath}
\usepackage{latexsym}
\numberwithin{equation}{section}

\def\R{{\bf R}}

\def\N{{\bf N}}

\def\d{\displaystyle}
\def\e{{\varepsilon}}
\def\o{\overline}
\def\wt{\widetilde}
\def\wh{\widehat}

\def\p{\partial}

\newtheorem{thm}{Theorem}[section]

\newtheorem{lem}{Lemma}[section]
\newtheorem{prop}{Proposition}[section]
\newtheorem{rem}{Remark}[section]

\title{Global existence for semilinear wave equations
with the critical blow-up term\\
in high dimensions
\footnote{This work is partially supported by
the Grant-in-Aid for Scientific Research (C) (No. 24540183),
Japan Society for the Promotion of Science.}
\vskip5pt
{\small }}
\author{
Hiroyuki Takamura
\footnote{Department of Complex and Intelligent Systems,
Faculty of Systems Information Science,
Future University Hakodate,
116-2 Kamedanakano-cho,
Hakodate, Hokkaido 041-8655, Japan.
e-mail: takamura@fun.ac.jp}
\quad
and
\quad
Kyouhei Wakasa
\footnote{The 2nd year of the doctor course,
Department of Mathematics,
Hokkaido University,
Sapporo 060-0810, Japan.
e-mail: wakasa@math.scihokudai.ac.jp}
}
\date{
\[
\begin{array}{ll}
\mbox{\footnotesize{\bf Keywords:}}
& \mbox{\footnotesize semilinear wave equation, high dimensions, critical exponent, lifespan}\\
\mbox{\footnotesize{\bf MSC2010:}}
& \mbox{\footnotesize primary 35L71, 35E15, secondary 35A01, 35A09, 35B33, 35B44}\\
\end{array}
\]
}

\begin{document}
\maketitle
\begin{abstract}
We are interested in almost global existence cases
in the general theory for nonlinear wave equations,
which are caused by critical exponents of nonlinear terms.
Such situations can be found in only three cases in the theory,
cubic terms in two space dimensions, quadratic terms in three space dimesions
and quadratic terms including a square of unknown functions itself
in four space dimensions.
Except for the last case, criterions to classify nonlinear terms
into the almost global, or global existence case,
are well-studied and known to be so-called null condition and non-positive condition.
\par
Our motivation of this work is to find such a kind of the criterion in four space dimensions.
In our previous paper, an example of the non-single term
for the almost global existence case is introduced. 
In this paper, we show an example of the global existence case.
These two examples have nonlinear integral terms
which are closely related to derivative loss due to high dimensions.
But it may help us to describe the final form of the criterion.
\end{abstract}

\section{Introduction}
\par
First we shall outline the general theory on the initial value problem
for fully nonlinear wave equations,
\begin{equation}
\label{GIVP}
\left\{
\begin{array}{l}
u_{tt}-\Delta u=H(u,Du,D_xDu) \quad \mbox{in}\quad\R^n\times[0,\infty),\\
u(x,0)=\e f(x),\ u_t(x,0)=\e g(x),
\end{array}
\right.
\end{equation}
where $u=u(x,t)$ is a scalar unknown function of space-time variables,
\[
\begin{array}{l}
Du=(u_{x_0},u_{x_1},\cdots,u_{x_n}),\ x_0=t,\\
D_xDu=(u_{x_ix_j},\ i,j=0,1,\cdots,n,\ i+j\ge1),
\end{array}
\]
$f,g\in C^\infty_0(\R^n)$ and $\e>0$ is a \lq\lq small" parameter.
We note that it is impossible to construct a general theory
for \lq\lq large" $\e$ due to blow-up results.
For example, see Glassey \cite{G73}, Levine \cite{Le74}, or Sideris \cite{Si84b}.
Let
\[
\wh{\lambda}=(\lambda;\ (\lambda_i),i=0,1,\cdots,n;
\ (\lambda_{ij}),i,j=0,1,\cdots,n,\ i+j\ge1). 
\]
Suppose that the nonlinear term $H=H(\wh{\lambda})$ is a sufficiently smooth function with
\[
H(\wh{\lambda})=O(|\wh{\lambda}|^{1+\alpha})
\]
in a neighborhood of $\wh{\lambda}=0$, where $\alpha\ge1$ is an integer.
Let us define the lifespan $\wt{T}(\e)$ of classical solutions of (\ref{GIVP}) by
\[
\begin{array}{ll}
\wt{T}(\e)=\sup\{t>0\ :\ \exists
&\mbox{a classical solution $u(x,t)$ of (\ref{GIVP})}\\
&\mbox{for arbitrarily fixed data, $(f,g)$.}\}.
\end{array}
\]
When $\wt{T}(\e)=\infty$, the problem (\ref{GIVP}) admits a global solution,
while we only have a local solution on $[0,\wt{T}(\e))$
when $\wt{T}(\e)<\infty$.
For local solutions, one can measure the long time stability
of a zero solution by orders of $\e$.
Because the uniqueness of the solution of (\ref{GIVP})
may yield that $\lim_{\e\rightarrow+0}\wt{T}(\e)=\infty$.
Such an uniqueness theorem can be found in Appendix of John \cite{J90} for example.
\par
In Chapter 2 of Li and Chen \cite{LC92},
we have long histories on the estimate for $\wt{T}(\e)$.
The lower bounds of $\wt{T}(\e)$ are summarized in the following table.
Let $a=a(\e)$ satisfy
\begin{equation}
\label{a}
a^2\e^2\log(a+1)=1
\end{equation}
and $c$ stands for a positive constant independent of $\e$.
Then,
due to the fact that it is impossible to obtain an $L^2$ estimate for $u$ itself
by standard energy methods, we have
\begin{center}
\begin{tabular}{|c||c|c|c|}
\hline
$\wt{T}(\e)\ge$  & $\alpha=1$ & $\alpha=2$ & $\alpha\ge3$\\
\hline
\hline
$n=2$ &
$\begin{array}{l}
ca(\e)\\
\quad\mbox{in general case},\\
c\e^{-1}\\
\quad\mbox{if}\ \d\int_{\R^2}g(x)dx=0,\\
c\e^{-2}\\
\quad\mbox{if}\ \partial^2_uH(0)=0
\end{array}
$
&
$\begin{array}{l}
c\e^{-6}\\
\quad\mbox{in general case},\\
c\e^{-18}\\
\quad\mbox{if}\ \partial^3_uH(0)=0,\\
\exp(c\e^{-2})\\
\quad\mbox{if}\ \partial^3_uH(0)=\partial^4_uH(0)=0\
\end{array}
$
& $\infty$ \\
\hline
$n=3$ &
$\begin{array}{l}
c\e^{-2}\\
\quad\mbox{in general case},\\
\exp(c\e^{-1})\\
\quad\mbox{if}\ \partial^2_uH(0)=0
\end{array}
$
& $\infty$ & $\infty$ \\\hline
$n=4$ &
$\begin{array}{l}
\exp(c\e^{-2})\\
\quad\mbox{in general case},\\
\infty\\
\quad\mbox{if}\ \partial^2_uH(0)=0
\end{array}$
& $\infty$ & $\infty$ \\
\hline
$n\ge5$ & $\infty$ & $\infty$ & $\infty$\\
\hline
\end{tabular} 
\end{center}
The result for $n=1$ is that
\begin{equation}
\label{lifespan_lower_onedim}
\wt{T}(\e)\ge
\left\{
\begin{array}{lllll}
\e^{-\alpha/2}\quad&\mbox{in general case},&\\
\d c\e^{-\alpha(1+\alpha)/(2+\alpha)}\quad&\mbox{if}\ \d \int_{\R}g(x)dx=0,&\\
c\e^{-\alpha} \quad &\mbox{if}\ \partial^\beta_uH(0)=0\ \mbox{for}\ 1+\alpha\le\forall\beta\le2\alpha.&
\end{array}
\right.
\end{equation}
For references on these results, see Li and Chen \cite{LC92}.
We shall skip to refer them here.
But we note that two parts in this table are different
from the one in Li and Chen \cite{LC92}.
One is the general case in $(n,\alpha)=(4,1)$.
In this part, the lower bound of $\wt{T}(\e)$ is
$\exp(c\e^{-1})$ in Li and Chen \cite{LC92}.
But later, it has been improved by Li and Zhou \cite{LZ95}.
Another is the case for $\partial^3_uH(0)=0$ in $(n,\alpha)=(2,2)$. 
This part is due to Katayama \cite{Kata01}.
But it is missing in Li and Chen \cite{LC92}.
Its reason is closely related to the sharpness of results in the general theory.
The sharpness is achieved by the fact that there is no possibility to improve
the lower bound of $\wt{T}(\e)$ in sense of order of $\e$
by blow-up results for special equations and special data.
It is expressed in the upper bound of $\wt{T}(\e)$
with the same order of $\e$ as in the lower bound.
On this matter, Li and Chen \cite{LC92} says that
all these lower bounds are known to be sharp except for $(n,\alpha)=(4,1)$.
But before this article,
Li \cite{Li91} says that $(n,\alpha)=(2,2)$ has also open sharpness
while the case for $\partial^3_uH(0)=0$ is still missing.
Li and Chen \cite{LC92} might have dropped the open sharpness in $(n,\alpha)=(2,2)$
by conjecture that $\partial^4_uH(0)=0$ is a technical condition.
No one disagrees with this observation
because the model case of $H=u^4$ has a global solution
in two space dimensions, $n=2$.
However, Zhou and Han \cite{ZH12} have obtained
this final sharpness in $(n,\alpha)=(2,2)$ by studying $H=u_t^2u+u^4$.
This puts Katayama's result into the table after 20 years from Li and Chen \cite{LC92}.
We note that Godin \cite{Go93} has showed the sharpness
of the case for $\partial^3_uH(0)=\partial^4_uH(0)=0$ in $(n,\alpha)=(2,2)$
by studying $H=u_t^3$.
This result has reproved by Zhou and Han \cite{ZH11_two}.
\par
We now turn back to another open sharpness of the general case in $(n,\alpha)=(4,1)$.
It has been obtained by our previous work, Takamura and Wakasa \cite{TW11},
by studying model case of $H=u^2$.
This part had been open more than 20 years
in the analysis on the critical case for model equations, $u_{tt}-\Delta u=|u|^p\ (p>1)$.
In this way, the general theory and its optimality have been completed. 
\par
After the completion of the general theory,
we are interested in the almost global existence,
namely, the case where $\wt{T}(\e)$ has an lower bound of the exponential function
of $\e$ with a negative power.
Such a case appears in $(n,\alpha)=(2,2),(3,1),(4,1)$
in the table of the general theory.
It is remarkable that Klainerman \cite{Kl86} and
Christodoulou \cite{Chri86} have independently found a special structure
on $H=H(Du,D_xDu)$ in $(n,\alpha)=(3,1)$
which guarantees the global existence.
This algebraic condition on nonlinear terms of derivatives
of the unknown function is so-called \lq\lq null condition".
It has been also established independently by Godin \cite{Go93} for $H=H(Du)$ 
and Katayama \cite{Kata95} for $H=H(Du,D_xDu)$ in $(n,\alpha)=(2,2)$.
The null condition had been supposed to be not sufficient
for the global existence in $(n,\alpha)=(2,2)$.
For this direction, Agemi \cite{A} proposed \lq\lq non-positive condition"
in this case for $H=H(Du)$.
This conjecture has been verified by Hoshiga \cite{Ho08} and Kubo \cite{K07} independently.
It might be necessary and sufficient condition to the global existence.
On the other hand, the situation in $(n,\alpha)=(4,1)$
is completely different from $(n,\alpha)=(2,2),(3,1)$
because $H$ has to include $u^2$.
\par
In our previous paper \cite{TW14}, we get the first attempt to clarify a criterion
on $H$ guaranteeing the global existence
by showing different blow-up example of $H$ from $u^2$ only.
More precisely, we have an almost global existence and its optimality
for an equation of the form
\begin{equation}
\label{integral_term1}
\begin{array}{ll}
u_{tt}-\Delta u=u^2
&\d-\frac{1}{\pi^2}\int_{0}^{t}d\tau
\int_{|\xi|\le 1}\frac{(u_tu)(x+(t-\tau)\xi,\tau)}{\sqrt{1-|\xi|^2}}d\xi\\
&\d-\frac{\e^2}{2\pi^2}\int_{|\xi|\le 1}\frac{f(x+t\xi)^2}{\sqrt{1-|\xi|^2}}d\xi
\end{array}
\end{equation}
in $\R^4\times[0,\infty)$.
We note that the third term in the right-hand side of (\ref{integral_term1})
can be neglected by simple modification.
One can say that this result is the first example of the blowing-up
of a classical solution to nonlinear wave equation with non-single 
and indefinitely signed term in high dimensions. 
We note that (\ref{integral_term1}) arises from a neglect
of derivative loss factors in Duhamel's term
for positive and single nonlinear term, $u^2$.
Therefore one can conclude that
derivative loss factors in Duhamel's term due to high dimensions
do not contribute to any order of $\e$ in the estimate of the lifespan.
\par
In this paper, we show that, in contrast with (\ref{integral_term1}),
another equation of the form
\begin{equation}
\label{integral_term2}
\begin{array}{ll}
u_{tt}-\Delta u=u^2
&\d-\frac{1}{2\pi^2}\int_0^t d\tau\int_{|\omega|=1}
(u_tu)(x+(t-\tau)\omega,\tau)dS_\omega\\
&\d-\frac{\e}{4\pi^2}\int_{|\omega|=1}(\e f^2+\Delta f+2\omega\cdot\nabla g)
(x+t\omega)dS_\omega
\end{array}
\end{equation}
admits a global classical solution in $\R^4\times[0,\infty)$.
Both the first integral terms in (\ref{integral_term1}) and (\ref{integral_term2})
look similar to each others.
The essential difference is that the second integral term in (\ref{integral_term2})
has linear terms of the initial data.
This part mainly comes from a neglect of derivative loss factors in the linear part.
Therefore one may say that
derivative loss factors in the linear part due to high dimensions
contribute to estimates of the lifespan.
\par
This paper is organized as follows.
In the next section, our main theorems are stated in more general situation
on space dimensions and nonlinear terms
as well as our motivation of this work by some inetgral equation.
In section 3, we investigate a relation between such an integral equation
and (\ref{integral_term2}).
The decay estimate of the linear part is studied in section 4.
The proof of the local existence appears in section 5,
for which {\it a priori} estimate in section 6 is required.
In the final section, blow-up result is proved
to show the optimality of the local existence.
\par
This work has begun since the second author
was in the 2nd year of the master course,
Graduate School of Systems Information Science,
Future University Hakodate.

\section{Main results}
This work is initiated by Agemi and Takamura \cite{AT94}
which attempts to make a new representation formula
of a solution of the following initial value problem
for inhomogeneous wave equations.
\begin{equation}
\label{LIVP}
\left\{
\begin{array}{l}
\partial_t^2u-\Delta u=F\quad  \mbox{in}\ \R^n\times[0,\infty),\\
u(x,0)=\e f(x),\ u_t(x,0)=\e g(x),\ x\in \R^n,
\end{array}
\right.
\end{equation}
where $u=u(x,t)$ is an unknown function,
$f,g$ and $F=F(x,t)$ are given smooth functions.
In \cite{AT94}, it has proved that, for $n\ge 3$, a smooth solution of (\ref{LIVP})
has to satisfy the following integral equation.
\begin{equation}
\label{AT_formula}
\begin{array}{ll}
(n-2)\omega_nu(x,t)=
&\d\e \int_{|\omega|=1}
\left\{t\omega\cdot \nabla f+(n-2)f+tg\right\}(x+t\omega)dS_{\omega}\\
&\d+(n-3)\int_{0}^{t}d\tau\int_{|\omega|=1}u_t(x+(t-\tau)\omega,\tau)dS_{\omega}\\
&\d+\int_{0}^{t}(t-\tau)d\tau\int_{|\omega|=1}F(x+(t-\tau)\omega,\tau)dS_{\omega},\\
\end{array}
\end{equation}
where $\omega_n$ is a measure of the unit sphere in $\R^n$, i.e. 
\[
\omega_{n}=\frac{2\pi^{n/2}}{\Gamma\left(n/2\right)}=
\left\{
\begin{array}{lc}
\d\frac{2(2\pi)^m}{(2m-1)!!} & \mbox{for}\ n=2m+1,\\
\d\frac{2\pi^{m+1}}{m!} & \mbox{for}\ n=2m+2,
\end{array}
\right.
(m=1,2,3,\cdots).
\]
In view of (\ref{AT_formula}),
neglecting the second term in the right-hand side,
we obtain a representation formula of a solution of some wave equation.
With a small modification, it may have the same initial data as in (\ref{LIVP}).
Our problem arises in this way. 
\par
In fact, let us define our integral equation of an unknown function $u$ by 
\begin{equation}
\label{AT_formula_new}
u(x,t)=\e V(x,t)+N(F)(x,t),
\end{equation}
where 
\begin{equation}
\label{linear_part}
V(x,t)=\frac{1}{\omega_n}\int_{|\omega|=1}
\left(\frac{t\omega\cdot\nabla f}{n-2}+f+tg\right)(x+t\omega)dS_\omega 
\end{equation}
and
\begin{equation}
\label{Duhamel_term}
N(F)(x,t)=\frac{1}{(n-2)\omega_n}
\int_{0}^{t}(t-\tau)d\tau \int_{|\omega|=1}F(x+(t-\tau)\omega,\tau)dS_\omega.
\end{equation}
Then, we have the following theorem.
\begin{thm}
\label{thm:AT_NW}
Let $n\ge3$. Assume that $f\in C^3(\R^n)$, $g\in C^2(\R^n)$
and $F\in C^2(\R^n\times [0,\infty))$.
Then, a solution of the integral equation (\ref{AT_formula_new})
satisfies the following initial value problem for inhomogeneous wave equation.
\begin{equation}
\label{AT_NWIVP}
\left\{
\begin{array}{ll}
\p_t^2u(x,t)-\Delta u(x,t)=F(x,t)-H(x,t) & \mbox{in}\quad\R^n\times[0,\infty),\\
u(x,0)=\e f(x),\ u_t(x,0)=\e g(x),& x\in \R^n,
\end{array}
\right.
\end{equation}
where $H$ is defined by
\begin{equation}
\label{AT_H}
\begin{array}{ll}
H(x,t)=&\d\frac{n-3}{(n-2)\omega_n}\int_{0}^{t}d\tau
\int_{|\omega|=1}(\partial_tF)(x+(t-\tau)\omega,\tau)dS_\omega\\
&\d+\frac{n-3}{(n-2)\omega_n}\int_{|\omega|=1}F(x+t\omega,0)dS_\omega\\
&\d+\frac{\e(n-3)}{(n-2)\omega_n}
\int_{|\omega|=1}\left\{\Delta f+(n-2)\omega\cdot\nabla g\right\}(x+t\omega)dS_\omega.
\end{array}
\end{equation}
\end{thm}
We shall make use of this theorem with $F(x,t)=u(x,t)^2$ and $n=4$.
The proof of this theorem appears in the next section.
\begin{rem}
The uniqueness of the solution of (\ref{AT_NWIVP})
with $F(x,t)=|u(x,t)|^p$ $(p\ge2)$ is open.
The restricted uniqueness theorem such as in Appendix 1 in John \cite{J90}
cannot be applicable because (99a) in \cite{J90} does not hold for this case.
\end{rem}
\begin{rem}
It is remarkable that Huygens' principle holds for $V$ in (\ref{linear_part})
even if the space dimension is even number.
See (\ref{Huygens_v}) below. 
Moreover, in view of (\ref{linear_part}) and (\ref{Duhamel_term}),
we need lower regularities on the data and inhomogeneous term
than those from $H\equiv0$ to obtain a classical solution.
\end{rem}
\par
In order to describe our main theorems,
let us define a lifespan $\wh{T}(\e)$ of the integral equation (\ref{AT_formula_new}) by
\[
\begin{array}{ll}
\wh{T}(\e)=\sup\{t>0\ :\ \exists & \mbox{a solution $u$ of (\ref{AT_formula_new})
with $F=F(u)$}\\
&\mbox{for arbitrarily fixed data, $(f,g)$.}\},
\end{array}
\]
where \lq\lq solution" means a classical solution of (\ref{AT_NWIVP}) for $p\ge2$,
or the $C^1$ solution of (\ref{AT_formula_new}) for $1<p<2$.
Our assumption on $F=F(s)$ is that
\begin{equation}
\left\{
\begin{array}{l}
\label{hypo_F_2}
\mbox{there exists a constant $A>0$ such that $F\in C^1(\R)$ satisfies}\\ 
\mbox{$|F^{(j)}(s)|\le A|s|^{p-j}\ (j=0,1)$ for $s\in\R,\ 1<p<2$},
\end{array}
\right.
\end{equation}
or
\begin{equation}
\left\{
\begin{array}{l}
\label{hypo_F_3}
\mbox{there exists a constant $A>0$ such that $F\in C^2(\R)$ satisfies}\\ 
\mbox{$|F^{(j)}(s)|\le A|s|^{p-j}$ $(j=0,1,2)$ for $s\in\R,\ p\ge2$}
\end{array}
\right.
\end{equation}
respectively. We also assume on the data that
\begin{equation}
\left\{
\begin{array}{l}
\label{hypo_data}
\mbox{at least one of $f\in C_0^{4}(\R^n)$ and $g\in C^{3}_{0}(\R^n)$ does not}\\
\mbox{vanish identically and have compact support}\\
\mbox{contained in $\{x\in\R^n\ :\ |x|\le k\}$ with some constant $k>1$.}
\end{array}
\right.
\end{equation}
We now introduce a critical number $p_1(n)$ as
a positive root of the following quadratic equation.
\begin{equation}
\label{p_1(n)}
\zeta(p,n)\equiv 2\left(1+(n-1)p-(n-2)p^2\right)=0.
\end{equation}
This is the analogy to Strauss' number defined by a positive root of
$\gamma(p,n)\equiv2+(n+1)p-(n-1)p^2=0$. See Remark \ref{rem:strauss} below.
\par
Then, we have the following lower bounds of the lifespan
which mean long time existences of the solution. 
\begin{thm}
\label{thm:main1}
Let $n\ge3$.
Assume that (\ref{hypo_F_2}), 
(\ref{hypo_F_3}) and (\ref{hypo_data}) are fulfilled.
Then there exists a positive constant $\e_0=\e_0(f,g,n,p,k)$ such that
the lifespan $\wh{T}(\e)$ satisfies
\begin{equation}
\label{lifespan_lower_main}
\begin{array}{ll}
\d \wh{T}(\e)=\infty
& \mbox{for}\ p>p_1(n),\\
\d \wh{T}(\e)\ge \exp\left(c\e^{-p(p-1)}\right)
& \mbox{for}\ p=p_1(n),\\
\d \wh{T}(\e)\ge c\e^{-2p(p-1)/\zeta(p,n)} 
& \mbox{for}\ 1<p<p_1(n)
\end{array}
\end{equation}
for any $\e$ with $0<\e \le\e_0$,
where $c$ is a positive constant independent of $\e$.
\end{thm}
\begin{rem}
\label{rem:strauss}
We note that
\[
p_1(n)=\frac{n-1+\sqrt{n^2+2n-7}}{2(n-2)}
\le p_0(n)=\frac{n+1+\sqrt{n^2+10n+7}}{2(n-1)}
\]
and that its equality holds if and only if $n=3$. 
Here $p_0(n)$ is Strauss' number on semilinear wave equations, $u_{tt}-\Delta u=|u|^p$.
See Strauss \cite{St81,St89} for this number.
Also see Takamura and Wakasa \cite{TW14} for references therein
on lifespan estimates for this equation.
Therefore the exponent $(n,p)=(4,2)$ is in the super critical case
for the equation (\ref{AT_NWIVP}) with $F(x,t)=|u(x,t)|^p$.
The key fact is that the linear part $V$ in (\ref{linear_part})
decays faster than that of a solution of the free wave equation. 
\end{rem}
\par
For the upper bounds of the lifespan,
our assumption on the data is the following.
\begin{equation}
\label{blowup_asm}
\left\{
\begin{array}{l}
\mbox{Let $f\equiv 0$, $g(x)=g(|x|)$ and $g\in C^{2}_{0}([0,\infty))$ satisfy that}\\
\mbox{(i) supp $g\subset\{x\in\R^n\ :\ |x|\le k\}$ with $k>0$},\\
\mbox{(ii) there exists $k_0$ such that $g(|x|)>0$ for $0<k_0<|x|<k$}.\\
\end{array}
\right.
\end{equation}
Then, we have the following theorem. 
\begin{thm}
\label{thm:main2}
Let $n\ge3$.
Assume that (\ref{blowup_asm}) is fulfilled.
Then there exists a positive constant $\e_0=\e_0(g,n,p,k)$ such that
the lifespan $\wh{T}(\e)$ satisfies
\begin{equation}
\label{lifespan_main}
\begin{array}{ll}
\d \wh{T}(\e)\le C\e^{-2p(p-1)/\zeta(p,n)}
\quad\mbox{for}\ 1<p<p_1(n),\\
\d \wh{T}(\e)\le\exp\left(C\e^{-p(p-1)}\right)
\quad\mbox{for}\ p=p_1(n)
\end{array}
\end{equation}
for any $\e$ with $0<\e \le\e_0$, where $C$ is a positive constant independent of $\e$.
\end{thm}
\par
The proofs of both Theorem \ref{thm:main1} and Theorem \ref{thm:main2}
are similar to those of our previous theorems in \cite{TW14}
which are based on John's iteration argument in a weighted $L^\infty$ space
by John \cite{J79}.
They are described after the next section.

\section{Proof of Theorem \ref{thm:AT_NW}}
First we shall prove the initial condition in (\ref{AT_NWIVP}).
It is trivial to get the first condition by setting $t=0$ in (\ref{AT_formula_new}).
Rewriting
\[
\omega\cdot\nabla f(x+t\omega)=\p_t(f(x+t\omega)),
\]
we have that
\begin{equation}
\label{AT_formula_new_diff_t}
\begin{array}{l}
u_t(x,t)\\
\d=\frac{\e}{\omega_n}\int_{|\omega|=1}\left\{\frac{((n-1)\p_t+t\p_t^2)f}{n-2}
+(1+t\p_t)g\right\}(x+t\omega)dS_{\omega}\\
\quad\d+\int_{0}^{t}d\tau\int_{|\omega|=1}
\frac{(1+(t-\tau)\p_t)F(x+(t-\tau)\omega,\tau)}{(n-2)\omega_n}dS_{\omega}.
\end{array}
\end{equation}
Therefore the second condition follows from setting $t=0$ in this equation.
\par
For the proof of the equation in (\ref{AT_NWIVP}),
we shall employ the well-known fact that a function $M(x,t)$ defined by 
\[
M(x,t)=\frac{1}{\omega_n}\int_{|\omega|=1}m(x+t\omega)dS_\omega
\]
for $m\in C^2(\R^n)$ satisfies the initial value problem of Darboux equation,
\begin{equation}
\label{Darboux}
\left\{
\begin{array}{l}
\d \left(\p_t^2-\Delta+\frac{n-1}{t}\p_t\right)M=0\\
M(x,0)=m(x),\ M_t(x,0)=0.
\end{array}
\right.
\end{equation}
Then it follows from (\ref{AT_formula_new_diff_t}) and (\ref{Darboux}) that 
\[
\begin{array}{ll}
u_{tt}(x,t)=
&\d\frac{\e}{\omega_n}\int_{|\omega|=1}
\left\{\frac{(1+t\p_t)\Delta f}{n-2}+(2\p_t+t\p^2_t)g\right\}
(x+t\omega)dS_\omega\\
&\d+\int_0^td\tau\int_{|\omega|=1}
\frac{(2\p_t+(t-\tau)\p^2_t)F(x+(t-\tau)\omega,\tau)}{(n-2)\omega_n}dS_\omega\\
&\d+\frac{F(x,t)}{n-2}.
\end{array}
\]
On the other hand, operating $\Delta$ to (\ref{AT_formula_new}) yields that
\[
\begin{array}{ll}
\Delta u(x,t)=
&\d\frac{\e}{\omega_n}\int_{|\omega|=1}
\left\{\left(1+\frac{t\p_t}{n-2}\right)\Delta f
+t\Delta g\right\}(x+t\omega)dS_\omega\\
&\d +\int_{0}^{t}(t-\tau)d\tau\int_{|\omega|=1}
\frac{\Delta F(x+(t-\tau)\omega,\tau)}{(n-2)\omega_n}dS_\omega.
\end{array}
\]
Therefore, it follows from (\ref{Darboux}) that 
\[
\begin{array}{l}
u_{tt}(x,t)-\Delta u(x,t)\\
\d=\frac{\e}{\omega_n}\int_{|\omega|=1}
\left\{\frac{3-n}{n-2}\Delta f+(2\p_t+t(\p_t^2-\Delta))g\right\}(x+t\omega)dS_\omega\\
\d+\int_0^td\tau\int_{|\omega|=1}
\frac{\{2\p_t+(t-\tau)(\p_t^2-\Delta)\}F(x+(t-\tau)\omega,\tau)}{(n-2)\omega_n}dS_{\omega}\\
\d+\frac{F(x,t)}{n-2}.
\end{array}
\]
Splitting $2\p_t$ into $(n-1)\p_t+(3-n)\p_t$
and making use of (\ref{Darboux}) again, we have that
\[
\begin{array}{l}
u_{tt}(x,t)-\Delta u(x,t)\\
\d=\frac{\e}{\omega_n}\int_{|\omega|=1}
\left\{\frac{3-n}{n-2}\Delta f+(3-n)(\p_t)g\right\}(x+t\omega)dS_\omega
+\frac{F(x,t)}{n-2}\\
\d+\frac{3-n}{(n-2)\omega_n}\int_0^td\tau\int_{|\omega|=1}
\p_t(F(x+(t-\tau)\omega,\tau))dS_\omega.
\end{array}
\]
Since
\[
\begin{array}{l}
\p_t(F(x+(t-\tau)\omega,\tau))\\
\d=(\p_tF)(x+(t-\tau)\omega,\tau)
-\p_{\tau}(F(x+(t-\tau)\omega,\tau))
\end{array}
\]
and
\[
\begin{array}{l}
\d\int_0^t\p_{\tau}
\left(\int_{|\omega|=1}F(x+(t-\tau)\omega,\tau)dS_\omega\right)d\tau\\
\d=\omega_nF(x,t)-\int_{|\omega|=1}F(x+t\omega,0)dS_\omega,
\end{array}
\]
we finally obtain that
\[
\begin{array}{l}
u_{tt}(x,t)-\Delta u(x,t)\\
\d=\d\frac{(3-n)\e}{(n-2)\omega_n}
\int_{|\omega|=1}\left\{\Delta f+(n-2)(\p_t)g\right\}(x+t\omega)dS_\omega\\
\d+\frac{3-n}{(n-2)\omega_n}\int_{|\omega|=1}F(x+t\omega,0)dS_\omega+F(x,t)\\
\d+\frac{3-n}{(n-2)\omega_n}\int_0^td\tau\int_{|\omega|=1}
(\p_tF)(x+(t-\tau)\omega,\tau)dS_\omega.
\end{array}
\]
This ends the proof of Theorem \ref{thm:AT_NW}. \hfill$\Box$
\section{Decay estimate of the linear part}
In this section, we get a space-time decay estimate of 
$V$ in (\ref{linear_part})
which plays an essential role to define our weighted $L^\infty$ space.
\begin{lem}
\label{lem:decay_est_v}
Under the same assumption as in Theorem \ref{thm:main1},
there exists a positive constant $C_{n,k}$
depending only on $n$ and $k$ such that
$V$ satisfies
\begin{equation}
\label{decay_est_v}
\begin{array}{l}
\d (t+|x|+2k)^{n-2}|\nabla_x^{\alpha}V(x,t)|\\
\d\le C_{n,k}\left(\sum_{|\beta|\le|\alpha|+2}\|\nabla_x^\beta f\|_{L^\infty(\R^n)}
+\sum_{|\gamma|\le|\alpha|+1}\|\nabla_x^\gamma g\|_{L^\infty(\R^n)}\right)
\end{array}
\end{equation}
for $|\alpha|\le2,\ (x,t)\in\R^n\times[0,\infty)$, and
\begin{equation}
\label{Huygens_v}
\mbox{supp}\ V\subset\{(x,t)\in\R^n\times[0,\infty)\ :\ -k\le t-|x|\le k\}.
\end{equation}
\end{lem}
\par\noindent
{\bf Proof.} First we note that the support property (\ref{Huygens_v})
immediately follows from the representation of $V$ in (\ref{linear_part}),
and that it is enough to prove the lemma for $|\alpha|=0$.
For (\ref{decay_est_v}) with $|\alpha|=0$,
one can employ the standard argument as in Lemma 3.2 in Agemi, Kubota and Takamura \cite{AKT94}.
\par
When $t\ge k$, taking into account of (\ref{Huygens_v}), one can make use of
\[
t^{n-1}\int_{|\omega|=1}|\varphi(x+t\omega)|dS_\omega\le\|\nabla_x\varphi\|_{L^1(\R^n)}
\quad\mbox{for}\ \varphi\in C_0^1(\R^n),\ t>0
\]
with
\[
t\ge\frac{1}{5}\left(t+|x|+2k\right).
\]
Hence we obtain that
\[
|V(x,t)|\le\frac{C_{n,k}}{(t+|x|+2k)^{n-2}}\left(\sum_{1\le|\beta|\le2}\|\nabla_x^\beta f\|_{L^1(\R^n)}
+\sum_{|\gamma|=1}\|\nabla_x^\gamma g\|_{L^1(\R^n)}\right)
\]
with some positive constant $C_{n,k}$ depending only on $n$ and $k$.
When $t\le k$, (\ref{linear_part}) yields that
\[
|V(x,t)|\le C_{n,k}\left(\sum_{|\beta|\le1}\|\nabla_x^\beta f\|_{L^\infty(\R^n)}
+\|g\|_{L^\infty(\R^n)}\right)
\]
with a different constant $C_{n,k}>0$.
Therefore the proof is completed. 
\hfill$\Box$
\section{Proof of Theorem \ref{thm:main1}}
Following Takamura and Wakasa \cite{TW14},
we prove Theorem \ref{thm:main1} in this section.
We note that its proof is similar to the one of odd dimensional case in \cite{TW14}
because of Huygens' principle for the linear part of the integral equation,
(\ref{Huygens_v}). 
It is obvious that the theorem follows from the following proposition.
\begin{prop}
\label{prop:lower_lifespan}
Let $n\ge3$.
Suppose that the assumptions (\ref{hypo_F_2}), (\ref{hypo_F_3}) and 
(\ref{hypo_data}) are fulfilled.
Then, there exists a positive constant $\e_0=\e_0(f,g,n,p,k)$ 
such that (\ref{AT_formula_new}) admits a unique solution $u\in C^1(\R^n\times[0,T])$
for $1<p<2$, $u\in C^2(\R^n\times[0,T])$ for $p\ge2$, as far as T satisfies
\begin{equation}
\label{lower_lifespan}
\begin{array}{ll}
\d T\le c\e^{-2p(p-1)/\zeta(p,n)} & \mbox{if}\ 1<p<p_1(n),\\
\d T\le\exp\left(c\e^{-p(p-1)}\right) & \mbox{if}\ p=p_1(n),\\
\d \mbox{there is no bound} & \mbox{if}\ p>p_1(n)
\end{array}
\end{equation}
for $0<\e \le \e_0$, where $c$ is a positive constant independent of $\e$.
\end{prop}
\par
The solution is constructed by almost the same way as in \cite{TW14}.
Actually, we shall set $U=u-\e V$ and rewrite (\ref{AT_formula_new})
with $F=F(u)$ into the following form.
\begin{equation}
\label{AT_formula_new'}
U=N(F(U+\e V)).
\end{equation}
Since $V$ exists globally in time,
we have to consider the lifespan of the solution of (\ref{AT_formula_new'}).
Let us define the sequence of functions, $\{U_m\}_{m\in\N}$ by
\[
U_m=N(F(U_{m-1}+U_0))\ \mbox{and}\ U_0=\e V.
\]
We also denote a weighted $L^\infty$ norm of $U$ by 
\[
\|U\|=\sup_{(x,t)\in\R^n\times[0,T]}\{w(|x|,t)|U(x,t)|\}
\]
with the weighted function
\[
\begin{array}{ll}
w(r,t)=
\left\{
\begin{array}{ll}
\tau_+(r,t)^{n-2}\tau_-(r,t)^{\o{q}}
& \d\mbox{if}\ p>\frac{n-1}{n-2},\\
\d\tau_+(r,t)^{n-2}\left(\log4\frac{\tau_+(r,t)}{\tau_-(r,t)}\right)^{-1}
& \d\mbox{if}\ p=\frac{n-1}{n-2},\\
\tau_+(r,t)^{n-2+\o{q}}
& \d\mbox{if}\ 1<p < \frac{n-1}{n-2},
\end{array}
\right.
\end{array}
\]
where we set
\[
\o{q}=(n-2)p-(n-1)
\]
and 
\[
\tau_+(r,t)=\frac{t+r+2k}{k},\quad \tau_-(r,t)=\frac{t-r+2k}{k}.
\]
\par\noindent
{\bf Proof of Proposition \ref{prop:lower_lifespan}.}
In view of Proposition 5.1 in \cite{TW14},
the proof of this proposition follows from the following {\it a priori} estimate.
\begin{lem}
\label{lem:apriori}
Let $n\ge3$ and $N$ be a linear integral operator defined in (\ref{Duhamel_term}).
Assume that $U,U_0\in C^0(\R^n\times[0,T])$ with
supp $U\subset\{(x,t)\in\R^n\times[0,T]\ :\ |x|\le t+k\}$,
supp $U_0\subset\{(x,t)\in\R^n\times[0,T]\ :\ t-k\le|x|\le t+k\}$,
and $\|U\|,\|\tau_+^{n-2}U_0w^{-1}\|<\infty$.
Then, there exists a positive constant $C_{n,\nu,p}$
depending on $n,\nu$ and $p$ such that 
\begin{equation}
\label{apriori}
\|N(|U_0|^{p-\nu}|U|^\nu)\|\le C_{n,\nu,p}k^2
\left\|\frac{\tau_+^{n-2}}{w}U_0\right\|^{p-\nu}\|U\|^\nu\o{E}_\nu(T) 
\end{equation}
for $0\le\nu\le p$, where $\o{E}_\nu$ is defined by 
\begin{equation}
\label{E_nu}
\o{E}_\nu(T)=
\left\{
\begin{array}{cl}
1 & \d\mbox{if}\ p>\frac{n-1}{n-2},\\
\d\left(\frac{2T+3k}{k}\right)^{\nu\delta} & \d\mbox{if}\ p=\frac{n-1}{n-2},\\
\d\left(\frac{2T+3k}{k}\right)^{-\nu\o{q}} 
& \d\mbox{if}\ 1<p<\frac{n-1}{n-2},
\end{array}
\right.
\end{equation}
for $0\le\nu<p$ with any $\delta>0$ and
\begin{equation}
\label{E_p}
\o{E}_p(T)=
\left\{
\begin{array}{cl}
1 & \mbox{if}\ p>p_1(n),\\
\d\log\frac{2T+3k}{k} & \mbox{if}\ p=p_1(n),\\
\d\left(\frac{2T+3k}{k}\right)^{\zeta(p,n)/2} 
& \mbox{if}\ \d 1<p<p_1(n).
\end{array}
\right.
\end{equation}
\end{lem}
This lemma is proved in the next section.
\par
The construction of the solution in our proposition is completely
same as in the proof of lower bounds of the lifespan
in odd space dimensions in the section 5 of Takamura and Wakasa \cite{TW14},
if $(n-1)/2$ in the exponent of $\tau_+$, $(n+1)/(n-1)$ in the definition of $E_\nu(T)$,
$q$, $p_0(n)$ and $\gamma(p,n)$ are
substituted by $(n-2)$, $(n-1)/(n-2)$, $\o{q}$, $p_1(n)$, and $\zeta(p,n)$
in all the questions respectively.
Therefore, Proposition \ref{prop:lower_lifespan}
immediately follows from Lemma \ref{lem:apriori}
which is proved in the next section.
\hfill$\Box$
\section{A priori estimates}
In this section we prove Lemma \ref{lem:apriori} which plays a key role 
in the proof of Theorem \ref{thm:main1}. The proof follows from the following 
basic estimate.
\begin{lem}{\rm\bf (Basic estimate)}
\label{lem:basicest}
Let $N$ be the linear integral operator defined by {\rm (\ref{Duhamel_term})} and 
$a_1\ge0$, $a_2 \in\R$ and $a_3\ge0$. Then, there exists a positive constant 
$C_{n,p,a_1,a_2,a_3}$ such that 
\begin{equation}
\label{basic_est}
\begin{array}{lll}
N\left\{\tau_{+}^{-(n-2)p+a_1}\tau_{-}^{a_2}
\left(\log(4\tau_{+}/\tau_{-})\right)^{a_3}\right\}(x,t)\\
\d \le C_{n,p,a_1,a_2,a_3}k^2w(r,t)^{-1}
\left(\frac{2T+3k}{k}\right)^{a_1}\o{E}_{a_1,a_2,a_3}(T)
\end{array}
\end{equation}
for $|x|\le t+k$, $t\in[0,T]$, where $\o{E}_{a_1,a_2,a_3}(T)$ 
is defined by 
\begin{equation}
\label{E_gen}
\o{E}_{a_1,a_2,a_3}(T)=
\left\{
\begin{array}{lll}
\d 1 & 
\mbox{if}\ \d a_2<-1\ \mbox{and}\ a_3=0,\\
\d \log\frac{2T+3k}{k} &
\mbox{if}\ \d a_2=-1\ \mbox{and}\ a_3=0,\\
\d \left(\frac{2T+3k}{k}\right)^{\delta a_3} &
\mbox{if}\ a_2\le-1\ \mbox{and}\ a_3>0,\\
\d \left(\frac{2T+3k}{k}\right)^{1+a_2} &
\mbox{if}\ \d a_2>-1,
\end{array}
\right.
\end{equation}
where $\delta$ stands for any positive constant.
\end{lem}
To prove this lemma, we shall employ the following lemma
which is established by fundamental identity for spherical means by John \cite{J55}.
\begin{lem}[John \cite{J55}]
\label{lm:Planewave}
Let $b\in C([0,\infty))$.
Then, the identity
\begin{equation}
\label{Planewave}
\begin{array}{ll}
\d \int_{|\omega|=1}b(|x+\rho \omega|)dS_\omega
\d = 2^{3-n}\omega_{n-1}(r\rho)^{2-n}\int_{|\rho-r|}^{\rho+r}\lambda b(\lambda)
h(\lambda,\rho,r)d\lambda,
\end{array} 
\end{equation}
holds for $x\in\R^n$, $r=|x|$ and $\rho>0$, where 
\begin{equation}
\label{h}
h(\lambda,\rho,r)
=\{\lambda^2-(\rho-r)^2\}^{(n-3)/2}\{(\rho+r)^2-\lambda^2\}^{(n-3)/2}.
\end{equation}
\end{lem}
For the proof of this lemma, see Lemma 4.1 in Agemi, Kubota and Takamura \cite{AKT94}.
In order to estimate $h(\lambda,\rho,r)$,
we shall make use of the following four inequalities.
\begin{lem}
\label{lm:h}
Let $h(\lambda,\rho,r)$ be the function defined by (\ref{h}).
Suppose that $|\rho-r|\le\lambda\le \rho+r$,
or equivalently $|\lambda-r|\le\rho\le\lambda+r$,
and $\rho\ge 0$.
Then the following inequalities hold.
\begin{eqnarray}
\label{h_1}
h(\lambda,\rho,r) &\le & 4^{n-3}r^{n-3}\lambda^{n-3},\\
\label{h_2}
h(\lambda,\rho,r) &\le & 2^{n-3}\rho^{n-3}r^{(n-3)/2}\lambda^{(n-3)/2},\\
\label{h_3}
h(\lambda,\rho,r) &\le & 8^{n-3}\rho^{n-3}r^{n-3},\\
\label{h_4}
h(\lambda,\rho,r) &\le & 2^{n-3}\rho^{n-3}\lambda^{n-3}.
\end{eqnarray}
\end{lem}
{\bf Proof.}
(\ref{h_1}), (\ref{h_2}) and (\ref{h_3}) are
due to Lemma 4.2 in Agemi, Kubota and Takamura \cite{AKT94}
with elementary computations.
(\ref{h_4}) is due to Lemma 2.2 in Georgiev \cite{Geo96}
with geometrical observation.
But one may prove (\ref{h_4}) also by elementally computation as follows.
\[
\begin{array}{l}
4\rho^2\lambda^2-\{\lambda^2-(\rho-r)^2\}\{(\rho+r)^2-\lambda^2\}\\
=\lambda^4+\{4\rho^2-(\rho+r)^2-(\rho-r)^2\}\lambda^2+(\rho-r)^2(\rho+r)^2\\
=(\lambda^2+\rho^2-r^2)^2\ge0.
\end{array}
\]
\hfill$\Box$\\
\par\noindent
{\bf Proof of Lemma \ref{lem:basicest}.} 
The proof is almost the same as the one in the estimates for $I_{odd}$ in Lemma 4.5 of Takamura and Wakasa \cite{TW14}. 
We denote various positive constants depending only on $n$ and $p$ by $C$ 
which may change at place to place.
By virtue of Lemma \ref{lm:Planewave}, we have that 
\[
N\left\{\tau_{+}^{-(n-2)p+a_1}\tau_{-}^{a_2}
\left(\log(4\tau_{+}/\tau_{-})\right)^{a_3}\right\}(x,t)=I(r,t),
\]
where we set 
\begin{equation}
\label{I}
\begin{array}{c}
\d I(r,t)=Cr^{2-n}\int_{0}^{t}(t-\tau)^{3-n}d\tau
\int_{|t-\tau-r|}^{t-\tau+r}\hspace{-20pt}
\tau_{+}(\lambda,\tau)^{-(n-2)p+a_1}\tau_{-}(\lambda,\tau)^{a_2}\times\\
\d\qquad\times\left(\log4\frac{\tau_{+}(\lambda,\tau)}{\tau_{-}(\lambda,\tau)}\right)^{a_3}
\lambda h(\lambda,t-\tau,r)d\lambda.
\end{array}
\end{equation}
We shall estimate $I(r,t)$ on three domains,
\[
\begin{array}{l}
D_1=\{(r,t)\ |\ r\ge t-r>-k\ \mbox{and}\ r\ge 2k\},\\
D_2=\{(r,t)\ |\ r\ge t-r>-k\ \mbox{and}\ r\le 2k\},\\
D_3=\{(r,t)\ |\ t-r\ge r\}.
\end{array}
\]
\par\noindent
(i) Estimate in $D_1$,
\par
Making use of (\ref{h_4}), we get 
\[
\begin{array}{l}
\d I(r,t)\d \le Cr^{2-n}\int_{0}^{t}d\tau
\int_{|t-\tau-r|}^{t+r-\tau}\lambda^{n-2}\times\\
\d\quad\times\tau_{+}(\lambda,\tau)^{-(n-2)p+a_1}\tau_{-}(\lambda,\tau)^{a_2}
\left(\log4\frac{\tau_{+}(\lambda,\tau)}{\tau_{-}(\lambda,\tau)}\right)^{a_3}d\lambda.
\end{array}
\]
Changing variables in the above integral by
\[
\alpha=\tau+\lambda,\ \beta=\tau-\lambda,
\]
we get
\[
\begin{array}{l}
\d I(r,t)\le Cr^{2-n}\int_{-k}^{t-r}\left(\frac{\beta+2k}{k}\right)^{a_2}d\beta
\int_{|t-r|}^{t+r}(\alpha-\beta)^{n-2}\times\\
\d\qquad\times\d\left(\frac{\alpha+2k}{k}\right)^{-(n-2)p+a_1}
\left(\log{4\frac{\alpha+2k}{\beta+2k}}\right)^{a_3}d\alpha.
\end{array}
\]
It follows from
\[
\frac{r}{k}=\frac{r+2r+r}{4k}\ge \frac{\tau_{+}(r,t)}{4}
\]
in $D_1$ that
\begin{equation}
\label{D_1est2}
\begin{array}{l}
\d I(r,t)\le C\tau_{+}(r,t)^{2-n}\left(\frac{t+r+2k}{k}\right)^{a_1}\int_{-k}^{t-r}\left(\frac{\beta+2k}{k}\right)^{a_2}d\beta\times\\
\d\qquad\times\int_{t-r}^{t+r}\d\left(\frac{\alpha+2k}{k}\right)^{-1-\o{q}}
\left(\log{4\frac{\alpha+2k}{\beta+2k}}\right)^{a_3}d\alpha.
\end{array}
\end{equation}
\par
When $a_3=0$, $\alpha$-integral in (\ref{D_1est2}) is dominated by 
\[
\left\{
\begin{array}{lll}
\d Ck\tau_{-}^{-\o{q}} & 
\mbox{if}\ \d p>\frac{n-1}{n-2},\\
\d k\log\frac{\tau_{+}}{\tau_{-}} &
\mbox{if}\ \d p=\frac{n-1}{n-2},\\
\d Ck\tau_{+}^{-\o{q}} &
\mbox{if}\ \d 1<p<\frac{n-1}{n-2}
\end{array}
\right.
\]
and $\beta$-integral in (\ref{D_1est2}) is dominated by 
\[
\left\{
\begin{array}{lll}
\d \frac{-k}{1+a_2} & 
\mbox{if}\ \d a_2<-1,\\
\d k\log\frac{t-r+2k}{k} &
\mbox{if}\ \d a_2=-1,\\
\d \frac{k}{1+a_2}\left(\frac{t-r+2k}{k}\right)^{1+a_2} &
\mbox{if}\ \d a_2>-1.
\end{array}
\right.
\]
(\ref{basic_est}) is now established for $a_3=0$.
\par
When $a_3>0$, we employ the following simple lemma. 
\begin{lem}
\label{lm:log}
Let $\delta>0$ be any given constant. Then, we have 
\[
\log X \le \frac{X^{\delta}}{\delta}\ for\ X\ge 1. 
\]
\end{lem}
The proof of this lemma follows from elementary computation.
We shall omit it.
Then, it follows from Lemma \ref{lm:log} that 
\[
\begin{array}{l}
\d I(r,t)\le C(4\delta^{-1})^{a_3}\tau_{+}(r,t)^{2-n}
\left(\frac{t+r+2k}{k}\right)^{a_1+\delta a_3}\times\\
\d\qquad\times\int_{-k}^{t-r}\left(\frac{\beta+2k}{k}\right)^{a_2-\delta a_3}d\beta
\d\int_{t-r}^{t+r}\left(\frac{\alpha+2k}{k}\right)^{-1-\o{q}}d\alpha.
\end{array}
\]
The $\alpha$-integral above can be estimated by the same manner in the case of $a_3=0$.
The $\beta$-integral is dominated by 
\begin{equation}
\label{beta_int_est_2}
\left\{
\begin{array}{lll}
\d \frac{-k}{1+a_2-\delta a_3} & 
\mbox{if}\ \d a_2\le-1,\\
\d \frac{k}{1+a_2-\delta a_3}\left(\frac{t-r+2k}{k}\right)^{1+a_2-\delta a_3} &
\mbox{if}\ \d a_2>-1
\end{array}
\right.
\end{equation}
with $\delta>0$ satisfying $1+a_2-\delta a_3>0$.
Therefore $I$ is bounded in $D_1$
by the quantity in the right-hand side of (\ref{basic_est}) as desired.
It is obvious that such a restriction on $\delta>0$ is finally removed from the statement.
\par\noindent
(ii) Estimate in $D_2$ or $D_3$.
\par
In this case, the proof is completely same as the one in the estimates for 
$I_{odd}$ in Lemma 4.5 in Takamura and Wakasa \cite{TW14}, 
if $(n-1)/2$ in the exponent of $\tau_{+}$ is substituted by $(n-2)$.
Because the key fact, $1-(n-2)p<0$, is also trivial. 
Therefore, the proof of Lemma \ref{lem:basicest} is now completed. 
\hfill$\Box$
\vskip10pt
\par\noindent
{\bf Proof of Lemma \ref{lem:apriori}.}
Due to Huygens' principle for the linear part $V$, (\ref{Huygens_v}),
one can replace $\tau_{-}$ by $\d \tau_{-}\chi_{\{-k\le t-r\le k\}}$ 
in (\ref{basic_est}) when $0\le \nu <p$.
Then, the integral with respect to the variable $\beta=\tau-\lambda$ is bounded.  
In order to establish  Lemma \ref{lem:apriori}, it is sufficient to show
\[
\left\{
\begin{array}{ll}
N\left(\tau_{+}^{-(n-2)(p-\nu)}w^{-\nu}\chi_{\{-k\le t-r\le k\}}\right)(x,t)
\le C_{n,\nu,p}k^2\o{E}_\nu(T)
&\mbox{for}\ 0\le \nu<p,\\
N(w^{-p})(x,t)\le C_{n,p,p}k^2\o{E}_p(T)
&\mbox{for}\ \nu=p.
\end{array}
\right.
\] 
\par
To this end, setting
\[
\left\{
\begin{array}{lll}
\d a_1=a_3=0,\ a_2=-\nu\o{q} & 
\mbox{if}\ \d p>\frac{n-1}{n-2},\\
\d a_1=a_2=0,\ a_3=\nu &
\mbox{if}\ \d p=\frac{n-1}{n-2},\\
\d a_1=-\nu \o{q},\ a_2=a_3=0 &
\mbox{if}\ \d 1<p<\frac{n-1}{n-2}
\end{array}
\right.
\]
for $0\le \nu<p$ and
\[
\left\{
\begin{array}{lll}
\d a_1=a_3=0,\ a_2=-p\o{q} & 
\mbox{if}\ \d p>\frac{n-1}{n-2},\\
\d a_1=a_2=0,\ a_3=p &
\mbox{if}\ \d p=\frac{n-1}{n-2},\\
\d a_1=-p\o{q},\ a_2=a_3=0 &
\mbox{if}\ \d 1<p<\frac{n-1}{n-2}
\end{array}
\right.
\]
for $\nu=p$ in (\ref{basic_est}), we have (\ref{apriori}). 
\hfill$\Box$
\section{Proof of Theorem \ref{thm:main2}}
In this section, we prove Theorem \ref{thm:main2}
which obviously follows from Proposition \ref{prop:lifespan_upper_odd} below. 
Its proof is almost the same as the one in odd dimensional case
of Theorem 2.2 in Takamura and Wakasa \cite{TW14}
once the similar iteration frame is established.
\begin{prop}
\label{prop:lifespan_upper_odd} 
Suppose that the assumptions of Theorem \ref{thm:main2} are fulfilled.
Let $u$ be a $C^0$-solution of (\ref{AT_formula_new}) in $\R^n\times[0,T]$.
Then, there exists a positive constant $\e_0=\e_0(g,n,p,k)$
such that $T$ cannot be taken as 
\begin{eqnarray}
\label{upper_lifespan_odd}
& T>\d \exp\left(c\e^{-p(p-1)}\right)& \mbox{if}\ p=p_1(n),\\
\label{upper_lifespan_sub_odd}
& T>\d c\e^{-2p(p-1)/\zeta(p,n)}& \mbox{if}\ 1<p<p_1(n) 
\end{eqnarray}
for $0<\e \le \e_0$, where $c$ is a positive constant independent of $\e$.
\end{prop}
\par\noindent
{\bf Proof.}
Similarly to the proof of Proposition 7.1 in \cite{TW14},
we may assume that the solution of (\ref{AT_formula_new}) is radially symmetric 
without loss of the generality.
Let $u=u(r,t)$ be a $C^0$-solution of
\begin{equation}
\label{rad_IE}
u=\e V+N(|u|^p)\quad\mbox{in}\quad(0,\infty)\times[0,T],
\end{equation}
where we set 
\begin{equation}
\label{V_rad}
V(r,t)=Cr^{2-n}t^{3-n}\int_{|t-r|}^{t+r}
\lambda g(\lambda)h(\lambda,t,r)d\lambda,
\end{equation}
\begin{equation}
\label{N_rad}
N(|u|^p)(r,t)\d=\o{C}r^{2-n}\int_{0}^{t}(t-\tau)^{3-n}d\tau
\int_{|t-\tau-r|}^{t-\tau+r}\hspace{-20pt}
\lambda h(\lambda,t-\tau,r)|u(\lambda,\tau)|^pd\lambda,
\end{equation}
where $C$ and $\o{C}$ are positive constants depending only on $n$. 
\vskip10pt
\par\noindent
{\bf [The 1st step] Inequality of $u$.}
\begin{lem}
\label{lem:lower_v}
Assume (\ref{blowup_asm}). 
Then there exists a positive constant $C_{n,g,k}>0$ such that for $t+k_0<r<t+k_1$
and $t\ge k_2$,
\begin{equation}
\label{lower_v}
V(r,t)\ge \frac{C_{n,g,k}}{r^{n-2}},
\end{equation}
where $\d k_1=\frac{k+k_0}{2}$ and $k_2=k-k_0$.
\end{lem}
{\bf Proof.} 
Let $t+k_0<r<t+k_1$ and $t\ge k_2/2$. Then, (\ref{V_rad}) gives us 
\[
V(r,t)\ge Cr^{2-n}t^{3-n}\int_{k_1}^{k}
\lambda g(\lambda)h(\lambda,t,r)d\lambda.
\]
Note that 
\[
\begin{array}{l}
r+t+\lambda\ge r,\ \lambda+r-t\ge \lambda,\\ 
r+t-\lambda\ge r+t-k\ge 2t+k_0-k\ge t,\ 
\lambda+t-r\ge \lambda-k_1 
\end{array}
\]
hold in the domain of the integral above for $t+k_0<r<t+k_1$ and $t\ge k_2$. 
Hence, we get 
\[
\begin{array}{ll}
V(r,t)&\d\ge Cr^{-(n-1)/2}t^{-(n-3)/2}
\int_{k_1}^{k}\lambda^{(n-1)/2}g(\lambda)(\lambda-k_1)^{(n-3)/2}d\lambda\\
&\d\ge C\left(\frac{k-k_1}{2}\right)^{(n-3)/2}r^{-(n-2)}
\int_{(k+k_1)/2}^{k}\lambda^{(n-1)/2}g(\lambda)d\lambda
\end{array}
\]
for $t+k_0<r<t+k_1$ and $t\ge k_2$.
Therefore we obtain (\ref{lower_v}).
\hfill$\Box$
\vskip10pt
\par
Making use of this estimate of $V$, we have the following iteration frame.
\begin{lem}
\label{lm:frame_odd}
Let $u$ be a $C^0$-solution of (\ref{rad_IE}).
Assume (\ref{blowup_asm}). 
Then $u$ in $\Sigma_0=\left\{(r,t):\ 2k\le t-r\le r\right\}$ satisfies
\begin{equation}
\label{frame1_odd}
\begin{array}{l}
\d u(r,t)\ge\frac{\o{C}2^{(n-3)/2}(t-r)^{(n-1)/2}}{r^{(3n-7)/2}}\times\\
\d\times\int\!\!\!\int_{R(r,t)}\{(t-r-\tau+\lambda)(t+r-\tau-\lambda)\}^{(n-3)/2}
|u(\lambda,\tau)|^pd\lambda d\tau+\\
\qquad\d+\frac{E_1(t-r)^{(3n-5)/2-(n-2)p}}{{r}^{(3n-7)/2}}\e^p,
\end{array}
\end{equation}
where $\o{C}$ is the one in (\ref{N_rad}),
\[
E_1=\frac{\o{C}C_{n,g,k}^p(k_1-k_0)}{(n-1)2^{(n-2)p-(3n-11)/2}}
\]
and
\[
R(r,t)=\left\{(\lambda,\tau):\ t-r\le\lambda,
\tau+\lambda\le t+r, 2k\le \tau-\lambda\le t-r\right\}.
\]
\end{lem}
\par\noindent
{\bf Proof.}
Comparing $L_{odd}$ in (4.7) of \cite{TW14} with radially symmetric form of $N$ in (\ref{V_rad}) of this paper,
the difference between the proof of Lemma 7.2 of \cite{TW14} and the one of this lemma
has to appear only in the second term, $I_2$, which arises from the estimate of the linear part.
In view of the proof of Lemma 7.2 in \cite{TW14},
the desired estimate immediately follows from simple replacement of $1-(n-1)p/2$ in the exponent of $\alpha-\beta$ by $1-(n-2)p$. 
\hfill$\Box$
\vskip10pt
\par\noindent
{\bf [The 2nd Step] Comparison argument.}
\par
Let us consider a solution $w$ of 
\begin{equation}
\label{v_equal}
\begin{array}{ll}
\d w(t-r)=
&\d\frac{\o{C}2^{(n-5)/2}(t-r)^{(n-1)/2}}{r^{(3n-7)/2}}
\int_{2k}^{t-r}(t-r-\beta)^{(n-3)/2}d\beta\\
&\d \times\int_{2(t-r)+\beta}^{t+r}(t+r-\alpha)^{(n-3)/2}|w(\beta)|^p d\alpha\\
&\d+\frac{E_1(t-r)^{(3n-5)/2-(n-2)p}}{2{r}^{(3n-7)/2}}\e^p.
\end{array}
\end{equation}
Then we have the following comparison lemma. 
\begin{lem}
\label{lem:comparison_arg}
Let $u$ be a solution of (\ref{rad_IE}) and $w$ be a solution of (\ref{v_equal}). 
Then, $u$ and $w$ satisfy 
\[
u>w\quad\mbox{in}\ \Sigma_0.
\]
\end{lem}
\par\noindent
{\bf Proof.}
Comparing the relation between $u$ in Lemma 7.3 of \cite{TW14} and $w$ in (7.6) of \cite{TW14}
with the one between $u$ in Lemma \ref{lm:frame_odd} and $w$ in (\ref{v_equal}),
one can find no difference in the structure of proofs of both Lemma 7.4 of \cite{TW14} and this lemma. 
\hfill$\Box$
\vskip10pt
\par
By definition of $w$ in (\ref{v_equal}), we have
\[
\begin{array}{ll}
w(\xi)\ge
&\d \frac{\o{C}\xi^{3-n}}{2^{n-1}}\int_{2k}^{\xi}
(\xi-\beta)^{(n-3)/2}|w(\beta)|^pd\beta\\
&\d \times\int_{2\xi+\beta}^{3\xi}(3\xi-\alpha)^{(n-3)/2}d\alpha
+\frac{E_1}{2^{(3n-5)/2}}\xi^{-\o{q}-(n-2)}\e^p\\
\end{array}
\]
in $\Gamma_0$, where we set
\[
\xi=\frac{r}{2},\ \Gamma_0=\{t-r=\xi, r\ge 4k\}.
\]
Hence we obtain that
\[
w(\xi)\ge\frac{\o{C}\xi^{3-n}}{2^{n-2}(n-1)}\int_{2k}^{\xi}
(\xi-\beta)^{n-2}|w(\beta)|^pd\beta+\frac{E_1\xi^{-\o{q}-(n-2)}}{2^{(3n-5)/2}}\e^p
\]
for $\xi\ge2k$.
Then, it follows from the setting
\[
W(\xi)=\xi^{\o{q}+n-2}w(\xi)
\] 
that
\begin{equation}
\label{frame2}
\d W(\xi)\ge D_n\xi^{\o{q}+1}\int_{2k}^{\xi}
\frac{(\xi-\beta)^{n-2}|W(\beta)|^pd\beta}{\beta^{(n-2)p+p\o{q}}}
+E_2\e^p\quad \mbox{for}\ \xi\ge 2k,
\end{equation}
where we set
\[
D_n=\frac{\o{C}}{2^{n-2}(n-1)},\ E_2=\frac{E_1}{2^{(3n-5)/2}}.
\]
\vskip10pt
\par\noindent
{\bf Iteration frame in the case of $p=p_1(n)$.}
\par
By virtue of (\ref{frame2}), we get
\begin{equation}
\label{frame_cri}
W(\xi)\ge D_n\int_{2k}^{\xi}\left(\frac{\xi-\beta}{\xi}\right)^{n-2}
\frac{|W(\beta)|^p}{\beta^{p\o{q}}}d\beta+E_2\e^p
\quad\mbox{for}\ \xi\ge 2k.
\end{equation}
The above inequality is the iteration frame for the critical case. 
This inequality is the same as the one in (7.8) in \cite{TW14}, 
if $q$ is substituted by $\o{q}$.
\vskip10pt
\par\noindent
{\bf Iteration frame in the case of $1<p<p_1(n)$.}
\par
Because of the fact that $-(n-2)p-p\o{q}<0$ for $n\ge3$, (\ref{frame2}) yields 
\begin{equation}
\label{frame_sub_odd}
\d W(\xi)\ge D_n\xi^{-(n-2)-p\o{q}}
\int_{2k}^{\xi}(\xi-\beta)^{n-2}|W(\beta)|^pd\beta
+E_2\e^p\quad \mbox{for}\ \xi\ge 2k.
\end{equation}
The above inequality is the iteration frame for the subcritical case. 
This inequality is the same as the one in (8,2) in \cite{TW14}, 
if $q$ is substituted by $\o{q}$.
\par
Making use of (\ref{frame_cri}) and (\ref{frame_sub_odd}),
one can obtain Proposition \ref{prop:lifespan_upper_odd}
immediately by the same argument in \cite{TW14}.
Therefore the proof of Theorem \ref{thm:main2} is now completed.
\hfill$\Box$
\bibliographystyle{plain}

\end{document}